\documentclass{llncs}

\usepackage{latexsym}
\usepackage{amsfonts}
\usepackage{graphicx}

\title{Thompson's group and public key cryptography}

\author{Vladimir Shpilrain\inst{1} \and Alexander Ushakov\inst{2}}
\institute{Department of Mathematics, The City  College of
New York, New York, NY 10031 \email{shpilrain@yahoo.com}
\thanks{Research of the first author was partially
supported by the NSF grant DMS-0405105.} \and Department of
Mathematics, CUNY Graduate Center, New York, NY 10016
\email{aushakov@mail.ru}}

\newtheorem{algorithm}{Algorithm}
\def\CR{{\mathcal R}}

\begin{document}

\maketitle

\begin{abstract}
Recently, several public key exchange protocols based on
symbolic computation in non-commutative (semi)groups were
proposed as a more efficient alternative to well
established protocols based on numeric computation.
Notably, the protocols due to Anshel-Anshel-Goldfeld and
Ko-Lee et al. exploited the {\it conjugacy search problem}
in groups, which is a  ramification of the discrete
logarithm problem. However, it is a prevalent opinion now
that the conjugacy search problem alone is unlikely to
provide sufficient level of security no matter what
particular group is chosen as a platform.

 In this paper we employ another problem (we call it the
{\it  decomposition problem}), which is more general than the
conjugacy search problem, and we suggest to use R. Thompson's group
as a platform. This group is well known in many areas of
mathematics, including algebra, geometry, and analysis. It also has
several properties that make it fit for cryptographic purposes.
 In particular, we show here that the word problem in Thompson's
group is solvable in almost linear time.
\end{abstract}

\section{Introduction}

One of the possible generalizations of the {\it discrete logarithm
problem} to arbitrary groups is the so-called
 {\it  conjugacy search  problem}: given two elements $a, b$ of a group $G$
 and the information that $a^x=b$
for some $x \in G$, find at least one particular element $x$ like
that. Here $a^x$ stands for $x^{-1}ax$. The (alleged) computational
difficulty of this  problem in some particular groups (namely, in
braid groups) has been  used in several group based cryptosystems,
most notably in \cite{AAG}  and \cite{KLCHKP}.
  It seems however now that the conjugacy search problem alone
  is unlikely to provide sufficient
level of security; see \cite{Shpilrain} and \cite{SU} for
explanations.

 In this paper we employ another problem, which generalizes
the conjugacy search problem, but at the same time
resembles the factorization problem which is at the heart
of the RSA cryptosystem. This problem which some authors
(see e.g. \cite{CKLHC}, \cite{KLCHKP}) call the {\it
decomposition problem} is as follows:

\medskip

 {\it Given an element $w$ of a (semi)group $G$, a subset $A \subseteq G$
 and an  element $x\cdot w\cdot y$, find elements $x', y' \in A$  such that
  $x'\cdot w\cdot y' = x\cdot w\cdot y$.}

\medskip

 The conjugacy search  problem (more precisely, its subgroup-restricted version
 used in \cite{KLCHKP}) is a special case of the
  decomposition problem  if one takes  $x=y^{-1}$.

 The usual factorization problem for integers used in the RSA cryptosystem
is also a special case of the  decomposition  problem  if one takes $w=1$ and $G={\bf Z}^\ast_p$,
the multiplicative (semi)group of integers modulo $p$.  It is therefore conceivable that with more
complex (semi)groups used as platforms, the corresponding cryptosystem may be more secure. At the
same time, in the group that we use in this paper (R. Thompson's group), computing (the normal
form of) a product of elements is faster than in ${\bf Z}^\ast_p$.

\medskip

 A key exchange protocol based on the general decomposition problem is
 quite straightforward (see e.g. \cite{KLCHKP}):
given two subsets $A, B \subseteq G$ such that $ab=ba$ for any $a
\in A, ~b \in B$, and given a public element  $w \in G$, Alice
selects private $a_1, a_2 \in A$ and sends the  element $a_1 w a_2$
to Bob. Similarly, Bob selects private $b_1, b_2 \in B$ and sends
the element $b_1 w b_2$ to Alice. Then Alice computes $K_A=a_1 b_1 w
b_2 a_2$, and Bob computes $K_B=b_1 a_1 w a_2 b_2$. Since
$a_ib_i=b_ia_i$ in $G$, one has  $K_A=K_B=K$ (as an element of $G$),
which is now Alice's and Bob's common secret key.

 In this paper, we suggest the following modification of this protocol
 which appears to be more secure (at least for our particular choice
of the platform) against so-called ``length based"
 attacks (see e.g. \cite{HS}, \cite{HT}), according to
our experiments (see our Section \ref{selection}). Given two subsets $A, B \subseteq G$ such that
$ab=ba$ for any $a \in A, ~b \in B$, and given a public element  $w \in G$, Alice selects private
$a_1 \in A$ and $b_1 \in B$ and sends the element $a_1 w b_1$ to Bob. Bob selects private $b_2 \in
B$ and $a_2  \in A$ and sends the element $b_2 w a_2$ to Alice. Then Alice computes $K_A=a_1 b_2 w
a_2 b_1$, and Bob computes $K_B=b_2 a_1 w b_1 a_2$. Since $a_ib_i=b_ia_i$ in $G$, one has
$K_A=K_B=K$ (as an element of $G$), which is now Alice's and Bob's common secret key.


 The group that we suggest to use as the platform for this protocol
 is Thompson's group $F$ well known in many areas of
mathematics, including algebra, geometry, and analysis. This group
is  infinite non-abelian. For us, it is important that Thompson's
group has the following nice presentation in terms of generators and
defining relations:

\begin{equation} \label{eq:inf_pres}
F = \langle x_0 , x_1 , x_2 , \ldots \mid x_i^{-1} x_k x_i =
x_{k+1} ~ (k>i) \rangle.
\end{equation}

 This presentation is infinite. There are also finite presentations
of this group; for example,

$$F = \langle x_0 , x_1 , x_2 , x_3, x_4
\mid x_i^{-1} x_k x_i = x_{k+1} ~ (k>i, ~k<4) \rangle,$$

\noindent but it is the infinite presentation above that allows for a
convenient normal form, so we are going to use that presentation
in our paper.

 For a survey on various properties of Thompson's group, we
refer to \cite{CFP}. Here we only give a description of the
``classical" normal form for elements of $F$.

The classical normal form for an element of Thompson's group is
a word of the form
\begin{equation} \label{eq:norm_form1}
x_{i_1} \ldots x_{i_s} x_{j_t}^{-1} \ldots x_{j_1}^{-1},
\end{equation}
such that the following two conditions are satisfied:
\begin{enumerate}
 \item[{\bf (NF1)}] $i_1 \le ...\le i_s$ and $j_1\le \ldots \le j_t$
 \item[{\bf (NF2)}] if both $x_i$ and $x_i^{-1}$ occur, then either
$x_{i+1}$ or $x_{i+1}^{-1}$  occurs, too.
\end{enumerate}

 We say that a word $w$ is in  {\em
seminormal form} if it is of the form (\ref{eq:norm_form1}) and
satisfies (NF1).

   We show in Section \ref{WP} that the time complexity of reducing a
word of length $n$ to the normal form in Thompson's group  is
$O(|n| \log |n|)$, i.e., is almost linear in $n$.

Another  advantage of cryptographic protocols based on symbolic
computation over those  based on computation with numbers is the
 possibility to generate a random word one symbol at a time.
For example, in RSA, one uses random prime numbers which obviously cannot be generated one digit
at a time but rather have to be precomputed, which limits the key space unless one wants to
sacrifice the efficiency. We discuss key generation in more detail in our Section \ref{selection}.
\medskip

\noindent {\it Acknowledgments.} We are grateful to V. Guba for helpful comments and to R.
Haralick for making a computer cluster  in his lab available for our computer experiments.

\section{The  protocol}
\label{protocol}

Let $F$ be Thompson's group given by its standard infinite
presentation (\ref{eq:inf_pres}) and $s \in \mathbb{N}$ a positive
integer. Define sets $A_s$ and $B_s$ as follows. The set $A_s$
consists of elements whose normal form is of the type
$$x_{i_1} \ldots x_{i_m} x_{j_m}^{-1} \ldots x_{j_1}^{-1},$$
i.e.  positive and negative parts are of the same length $m$, and
\begin{equation}\label{eq:condition}
i_k-k < s \mbox{ and } j_k-k < s \mbox{ for every } k = 1,\ldots,s.
\end{equation}
The set $B_s$ consists of elements represented by words in
generators $x_{s+1},x_{s+2},\ldots$. Obviously, $B_s$ is a subgroup
of $F$.

\begin{proposition}
Let $a \in A_s$ and $b\in B_s$. Then $ab = ba$ in the group $F$.
\end{proposition}

\begin{proof}
Let $a = x_{i_1} \ldots x_{i_m} x_{j_m}^{-1} \ldots x_{j_1}^{-1}$
and $b = x_{k_1}^{\varepsilon_{1}} \ldots x_{k_l}^{\varepsilon_{l}}$
where $k_q > s$ for every $q=1,\ldots,l$. By induction on $l$ and
$m$ it is easy to show that in the group $F$ one has

$$a b = b a = x_{i_1} \ldots x_{i_m} \delta_m(b)  x_{j_m}^{-1} \ldots x_{j_1}^{-1},$$
where $\delta_M$ is the operator that increases indices of all
generators by $M$ (see also our Section \ref{WP}). This
establishes the claim.
\end{proof}

\begin{proposition}
Let $s \ge 2$ be an integer. The set   $A_s$ is a subgroup of $F$
generated by $x_0 x_1^{-1}, \ldots, x_0 x_s^{-1}.$
\end{proposition}

\begin{proof}
The set $A_s$ contains the identity and is clearly closed under
taking inversions, i.e., $A_s = A_s^{-1}$. To show that $A_s$ is
closed under multiplication we take two arbitrary normal forms from
$A_s$:
$$u = x_{i_1} \ldots x_{i_m} x_{j_m}^{-1} \ldots x_{j_1}^{-1}$$
and
$$v = x_{p_1} \ldots x_{p_l} x_{q_l}^{-1} \ldots x_{q_1}^{-1}$$
and show that the normal form of $uv$ belongs to $A_s$. First, note
that since the numbers of   positive and negative letters in $uv$
are equal,
 the lengths of the positive and negative parts in the normal
form of $uv$ will be equal, too (see the rewriting system in the beginning of our Section
\ref{WP}). Thus, it remains to show that the property (\ref{eq:condition}) of indices in the
normal form of $u v$ is satisfied. Below we sketch the proof of this claim.

Consider the subword in the middle of the product $uv$ marked below:
$$uv = x_{i_1} \ldots x_{i_m} \left( x_{j_m}^{-1} \ldots x_{j_1}^{-1}    x_{p_1} \ldots x_{p_l} \right) x_{q_l}^{-1} \ldots x_{q_1}^{-1}$$
and find a seminormal form for it using relations of $F$ (move
positive letters to the left and  negative letters to the right
starting in the middle of the subword). We refer the reader to
Algorithm \ref{al:semi_norm_form_NP} in Section \ref{WP} for more
information on how this can be done. Denote the obtained word by
$w$. The word $w$ is the  product of a positive and a negative word:
$w = p n$. By induction on $l+m$ one can show that both $p$ and $n$
  satisfy the condition (\ref{eq:condition}).

Then we find normal forms for words $p$ and $n$ using   relations of
$F$ (for $p$ move   letters with smaller indices to the left of
letters with bigger indices, and for $n$ move letters with smaller
indices  to the right of   letters with bigger indices). By
induction on the number of operations thus performed, one can show
that the obtained words $p'$ and $n'$ satisfy the condition
(\ref{eq:condition}). Therefore, the word $w' = p' n'$ is a
seminormal form of $uv$ satisfying the condition
(\ref{eq:condition}).

Finally, we remove those  pairs of generators in $w'$ that
contradict the property (NF2) (we refer the reader to our Algorithm
\ref{al:semi_norm_form_erase} for more information). Again, by
induction on the number of ``bad pairs", one can show that the
result will satisfy the condition (\ref{eq:condition}). Therefore,
$uv$ belongs to $A_s$, i.e., $A_s$ is closed under multiplication,
  and therefore, $A_s$ is a subgroup.

Now we show that the set of words $\{x_0 x_1^{-1}, \ldots, x_0
x_s^{-1}\}$ generates the subgroup $A_s$. Elements   $\{x_0
x_1^{-1}, \ldots, x_0 x_s^{-1}\}$ clearly belong to $A_s$. To show
the inclusion $A_s \le \langle x_0 x_1^{-1}, \ldots, x_0 x_s^{-1}
\rangle$, we construct the Schreier graph of $\langle x_0 x_1^{-1},
\ldots, x_0 x_s^{-1} \rangle$ (depicted in Figure
\ref{fi:schreier_graph}) and see that any word from $A_s$ belongs to
the subgroup on the right.
\begin{figure}[htbp]
\centerline{ \includegraphics[scale=0.6]{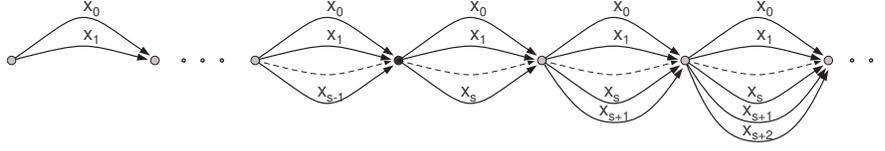} }
\caption{\label{fi:schreier_graph} The Schreier graph of the
subgroup $H = \langle x_0 x_1^{-1}, \ldots, x_0 x_s^{-1} \rangle$.
The black dot denotes the right coset corresponding to $H$.}
\end{figure}

\end{proof}

Now we give  a formal description of the  protocol based on the
decomposition problem mentioned in the Introduction.

\medskip

\noindent  {\bf (0)} Fix two positive integers $s, M$ and a word $w
= w(x_0,x_1,\ldots)$.
\medskip

\noindent  {\bf (1)} Alice randomly selects private elements $a_1
\in A_s$ and $b_1 \in B_s$. Then she reduces the  element $a_1 w
b_1$ to the normal form and sends the result to Bob.

\medskip

\noindent  {\bf (2)} Bob randomly selects private elements  $b_2 \in
B_s$ and $a_2 \in A_s$. Then he reduces the element $b_2 w a_2$  to
the normal form and sends the result to Alice.

\medskip

\noindent  {\bf (3)}  Alice computes $K_A= a_1 b_2 w a_2 b_1 = b_2 a_1 w b_1 a_2$, and Bob
computes $K_B= b_2 a_1 w b_1 a_2$. Since $a_ib_i=b_ia_i$ in $F$, one has $K_A=K_B=K$ (as an
element of $F$), which is now Alice's and Bob's common secret key.

\section{Parameters and key generation}
\label{selection}
\medskip

In practical key exchange we suggest to choose the following
parameters.

\medskip\noindent{\bf (1)} Select (randomly and uniformly) the parameter $s$
from the interval $[3,8]$ and the parameter $M$ from the set
$\{256,258,\ldots,318,320\}$.

\medskip\noindent{\bf (2)} Select the ``base" word $w$ as a
product of generators
$$S_W = \{x_{0},x_{1},\ldots,x_{s+2}\}$$
and their inverses. This is done the following way. We start with
the empty word $v_0$. When we have a current word $v_i$, we multiply
it on the right by a generator from $S_B^{\pm 1}$ and compute the
normal form of the product. The obtained word is denoted by
$v_{i+1}$. We continue this process until the obtained word
$v_{i+1}$ has length $M$.

\medskip\noindent{\bf (3)} Select   $a_1$ and $a_2$ as products of
words from
$$S_A = \{ x_0 x_1^{-1}, \ldots, x_0 x_s^{-1} \}$$
and their inverses. This  is done essentially the same way as above for $w$. We start with the
empty  word $u_0$. Let $u_i$ be the currently constructed word of length less than $M$. We
multiply $u_i$ on the right by a randomly chosen word from $S_A^{\pm 1}$ and compute the normal
form of the product. Denote the obtained normal form by $u_{i+1}$. Continue this process until the
obtained word $u_{i+1}$ has length $M$.

\medskip\noindent{\bf (4)} Select   $b_1$ and $b_2$ as
products of generators from
$$S_B = \{x_{s+1},x_{s+2},\ldots,x_{2s}\}$$
and their inverses. To  do that, start with the empty word $v_0$.
Multiply a current word $v_i$  on the right by a generator from
$S_B^{\pm 1}$ and compute the normal form of the product. Denote the
obtained  word  by $v_{i+1}$.  Continue this process until the
obtained word $v_{i+1}$ has length $M$.
\medskip

We would like to point out that the key space in the proposed scheme is exponential in   $M$; it
is easy to see that $|A_s(M)| \ge \sqrt{2}^M$.

The parameters above were chosen in such a way to prevent a  length-based attack.
 Note that for Thompson's  group,  a length-based attack could be a threat
since the normal form of any element represents a geodesic in the Cayley graph of $F$. Since ideas
behind length-based attacks were never fully described, we present below  a typical algorithm
(adapted to our situation) implementing such  an   attack (Algorithm \ref{al:length_attack}).

Define  a directed labelled graph $\Gamma = (V(\Gamma),E(\Gamma))$ as follows:
\begin{itemize}
 \item The set of vertices $V(\Gamma)$ corresponds to the set of all  elements of the group
 $F$.

\item The set of edges $E(\Gamma)$ contains edges $v_1
\stackrel{(w_1,w_2)}{\longrightarrow} v_2$ such that $v_2 = w_1 v_1 w_2$ in the group $F$, with
 labels of two types:
\begin{itemize}
 \item $(w_1,1)$, where $w_1 \in S_A^{\pm 1}$.
 \item $(1,w_2)$, where $w_2 \in S_B^{\pm 1}$.
\end{itemize}
\end{itemize}
For an element $w \in F$ denote by $\Gamma_w$ the connected component of $\Gamma$ containing $w$.
From the description of the protocol it follows that $w$ and the element $w' = a_1 w b_1$
transmitted by Alice to Bob belong to $\Gamma_w = \Gamma_{w'}$, and breaking Alice's key is
equivalent to finding a label of a path from $w$ to $w'$ in $\Gamma_w$.

To test our protocol, we performed a series of experiments. We randomly generated keys (as
described above) and ran  Algorithm \ref{al:length_attack} (see below) on them. Algorithm
\ref{al:length_attack} keeps constructing $\Gamma_w$ and $\Gamma_{w'}$ until a shared element is
found. The sets $S_w$ and $S_{w'}$ in the algorithm accumulate   constructed parts of the graphs
$\Gamma_{w}$ and $\Gamma_{w'}$. The sets $M_w \subset S_w$ and $M_{w'} \subseteq S_{w'}$ are
called the sets of marked vertices and are used to specify  vertices that are worked out.
\begin{algorithm} {\em (Length-based attack)} \label{al:length_attack}\\
{\sc Input.} The original public  word $w$ and the  word $w'$ transmitted by Alice.\\
{\sc Output.} A pair of words $x_1\in S_A$, $x_2\in S_B$
such that $w' = x_1 w x_2$.\\
{\sc Initialization.} Put $S_w = \{w\}$, $S_{w'} = \{w'\}$, $M_w = \emptyset$, $M_{w'} =
\emptyset$. \\
{\sc Computations.}
\begin{enumerate}
\item[A.] Find a shortest word $u \in S_w \setminus M_w$.

\item[B.] Multiply $u$ by elements $S_A^{\pm 1}$ on the
left and by elements $S_B^{\pm 1}$ on the right and add each result into $S_w$ with the edges
labelled accordingly.

\item[C.] Add $u$ into $M_w$.

\item[D.] Perform the steps A--C with $S_{w}$ and  $M_{w}$ replaced by $S_{w'}$ and  $M_{w'}$, respectively.

\item[E.] If $S_w \cap S_{w'} = \emptyset$ then goto A.

\item[F.] If there is $\overline{w} \in S_w \cap
S_{w'}$ then find a path in $S_w$ from $w$ to
$\overline{w}$ and a path in $S_{w'}$ from $\overline{w}$
to $w'$. Concatenate them and output the label of the
result.

\end{enumerate}
\end{algorithm}

We performed a series of   tests implementing  this length-based attack; in each test we let the
program to run overnight.   None of the programs gave a result, i.e., the success rate of the
length-based attack in our tests was 0.

\section{The word problem in Thompson's group}
\label{WP}
\medskip

 In this section, we show that the time complexity of reducing a
word of length $n$ to the normal form in Thompson's group $F$  is
$O(|n| \log |n|)$, i.e., is almost linear in $n$. Our algorithm is
in two independent parts: first we reduce a given word to a
seminormal form (Algorithm \ref{al:semi_norm_form2}), and then
further reduce it to the normal form by eliminating ``bad pairs"
(Algorithm \ref{al:semi_norm_form_erase}). We also note that crucial
for Algorithm \ref{al:semi_norm_form2} is Algorithm
\ref{al:semi_norm_form}  which computes a seminormal form of a
product of two seminormal forms. Our strategy for computing a
seminormal form of a given $w \in F$ is therefore recursive
(``divide and conquer"): we split the word $w$ into two halves:
$w=w_1w_2$, then compute  seminormal forms of $w_1$ and $w_2$, and
then use Algorithm \ref{al:semi_norm_form} to compute a seminormal
form of $w$.

 Recall that Thompson's group $F$ has the following infinite presentation:
$$F = \langle x_0 , x_1 , x_2 , \ldots \mid x_i^{-1} x_k x_i =
x_{k+1} ~ (k>i) \rangle.$$

The classical normal form for an element of Thompson's group (see
\cite{CFP} for more information) is described in the Introduction.

Let us denote by $\rho(w)$ the normal form for $w\in F$; it is
unique for a given element of $F$. Recall that we say that a word
$w$ is in {\em seminormal form} if it is of the form
(\ref{eq:norm_form1}) and satisfies (NF1) (see the Introduction). A
seminormal form is not unique. As usual, for a word $w$ in the
alphabet $X$ by $\overline{w}$ we denote the corresponding freely
reduced word.

As mentioned above, the normal form for an element of Thompson's
group can be computed in two steps:
\begin{enumerate}
 \item[1)] Computation of a seminormal form.
 \item[2)] Removing ``bad pairs", i.e., pairs $(x_i,x_i^{-1})$ for which the
 property (NF2) fails.
\end{enumerate}
The first part is achieved (Lemma
\ref{le:R_term}) by using the following rewriting system (for all
pairs $(i,k)$ such that $i<k$):
$$
\begin{array}{lll}
x_k x_i               & \rightarrow & x_i x_{k+1}           \\
x_k^{-1} x_i          & \rightarrow & x_i x_{k+1}^{-1}      \\
x_i^{-1} x_k          & \rightarrow & x_{k+1} x_i^{-1}      \\
x_i^{-1} x_k^{-1}     & \rightarrow & x_{k+1}^{-1} x_i^{-1} \\
\end{array}
$$
and, additionally, for all $i\in \mathbb{N}$
$$
\begin{array}{lll}
x_i^{-1} x_i & \rightarrow & 1 \\
\end{array}
$$
We denote this system of rules by $\CR$. It is straightforward to
check (using the confluence test, see \cite[Proposition 3.1]{Sims})
that $\CR$ is confluent. The following lemma is obvious.

\begin{lemma} \label{le:R_term}
$\CR$ terminates with a seminormal form. Moreover, a word is in a
seminormal form if and only if it is $\CR$-reduced.
\end{lemma}

Let us now examine the action of $\CR$ more closely. This action  is
  similar to  sorting a list of numbers,
but with two  differences: indices of generators may increase, and
some generators may disappear.

By Lemma \ref{le:R_term}, for any word $w$ in generators of $F$, the
final result of rewrites by $\CR$ is a seminormal form. Therefore,
to compute a seminormal form we implement rewrites by $\CR$. We do
it in a special manner in Algorithm \ref{al:semi_norm_form} in order
to provide the best performance. For convenience we introduce a
parametric function $\delta_\varepsilon, ~\varepsilon \in
\mathbb{Z},$  defined on the set of all words in the alphabet
$\{x_0^{\pm 1},x_1^{\pm 1},\ldots\}$ by
$$x_i^{\pm 1}  \stackrel{\delta_\varepsilon}{\mapsto} x_{i+\varepsilon}^{\pm 1}.$$

The function $\delta_\varepsilon$ may not be defined for some
negative $\varepsilon$ on a given  word  $w=w(x_{i_1}^{\pm 1},
x_{i_2}^{\pm 1}, \ldots )$,
 but when it is used, it is assumed that the function is defined.

\subsection{Merging seminormal forms}

Consider first the case where a word $w$ is a product of $w_1$ and
$w_2$ given in  seminormal forms. Let $w_1 = p_1 n_1$ and $w_2 = p_2
n_2$, where $p_i$ and $n_i$ ($i=1,2$) are the positive and negative
parts of $w_i$. Clearly, one can arrange the rewriting process for
$p_1 n_1 p_2 n_2$ by $\CR$ the following way:
\begin{enumerate}
 \item[1)] Rewrite the subword $n_1 p_2$ of $w$ to a seminormal form $p_2'
 n_1'$. Denote by $w' = p_1 p_2' n_1' n_2$ the obtained result.
 \item[2)] Rewrite the positive subword $p_1 p_2'$ of $w'$ to a seminormal form $p$.
 Denote by $w'' = p n_1' n_2$ the obtained result.
 \item[3)] Rewrite the negative  subword $n_1' n_2$ of $w''$ to a seminormal form $n$.
 Denote by $w''' = p n$ the obtained result.
\end{enumerate}
The word $w''' = pn$ is clearly in a seminormal form and $w =_F
w'''$.  This process can be depicted as follows:
$$
\begin{array}{c}
p_1 \underbrace{n_1 p_2} n_2 \\
\Downarrow \\
\underbrace{p_1 p_2'} \underbrace{n_1' n_2} \\
\Downarrow \\
p n
\end{array}
$$

The next algorithm performs the first rewriting step from the scheme
above, and the following  Lemma \ref{le:merging_correct_NP} asserts
that it  correctly performs the first step in linear time.

\begin{algorithm} \label{al:semi_norm_form_NP}{\em (Seminormal form of a product of negative and positive seminormal
forms)}\\
{\sc Signature.} $w = Merge_{-,+}(n,p , \varepsilon_1 , \varepsilon_2 )$.\\
{\sc Input.} Seminormal forms $n$ and $p$ (where $n = x_{j_t}^{-1}
\ldots x_{j_1}^{-1}$ and $p = x_{i_1} \ldots x_{i_s}$),
and numbers $\varepsilon_1, \varepsilon_2 \in \mathbb{Z}$.\\
{\sc Output.} Seminormal form $w$ such that $w =_F \delta_{\varepsilon_1}(n) \delta_{\varepsilon_2}(p)$.\\
{\sc Computations.}

\begin{enumerate}

 \item[A)] If $s=0$ or $t=0$ then output a product $n p$.

 \item[B)] If $j_1+\varepsilon_1 = i_1+\varepsilon_2$ then
\begin{enumerate}
 \item[1)] Compute $w = Merge_{-,+}(x_{j_t}^{-1} \ldots x_{j_2}^{-1},x_{i_2} \ldots
 x_{i_s} , \varepsilon_1 , \varepsilon_2)$.
 \item[2)] Output $w$.
\end{enumerate}

 \item[C)] If $j_1+\varepsilon_1 < i_1+\varepsilon_2$ then
\begin{enumerate}
 \item[1)] Compute $w = Merge_{-,+}(x_{j_t}^{-1} \ldots x_{j_2}^{-1},x_{i_1} \ldots
 x_{i_s} , \varepsilon_1 , \varepsilon_2+1 )$.
 \item[2)] Output $w x_{j_1+\varepsilon_1}^{-1}$.
\end{enumerate}

 \item[D)] If $j_1+\varepsilon_1 > i_1+\varepsilon_2$ then
\begin{enumerate}
 \item[1)] Compute $w = Merge_{-,+}(x_{j_t}^{-1} \ldots x_{j_1}^{-1},x_{i_2} \ldots
 x_{i_s} , \varepsilon_1+1 , \varepsilon_2 )$.
 \item[2)] Output $x_{i_1+\varepsilon_2} w$.
\end{enumerate}

\end{enumerate}
\end{algorithm}

\begin{lemma} \label{le:merging_correct_NP}
For any seminormal forms $n$ and $p$ (where $n = x_{j_t}^{-1} \ldots x_{j_1}^{-1}$ and $p =
x_{i_1} \ldots x_{i_s}$) and numbers $\varepsilon_1, \varepsilon_2 \in \mathbb{Z}$ the output $w =
Merge_{-,+}(n,p,  , \varepsilon_1 , \varepsilon_2 )$ of Algorithm \ref{al:semi_norm_form_NP} is a
seminormal form for $\delta_{\varepsilon_1}(n) \delta_{\varepsilon_2}(p)$. Furthermore,
  the time complexity required to compute $w$ is bounded
  by $C(|n|+|p|)$ for some constant $C$.
\end{lemma}

\begin{proof}

Since in each iteration we perform the constant number of elementary
steps and in each subsequent iteration the sum $|n|+|p|$ is
decreased by one, the time complexity of Algorithm
\ref{al:semi_norm_form_NP} is linear.

We prove correctness of Algorithm \ref{al:semi_norm_form_NP} by
induction on $|n|+|p|$. Assume that $|n|+|p| = 0$. Then at step A)
we get output $w = n p$ which is an empty word. Clearly, such $w$ is
a seminormal form for $n p$, so the base of  induction is done.

Assume that $|n|+|p| = N+1$ and for any shorter word the statement
is true. Consider four cases. If $|n| = 0$ or $|p| = 0$ then one of
the words is trivial and, obviously, the product $n p$ is a correct
output for this case. If $j_1+\varepsilon_1 = i_1+\varepsilon_2$
then $x_{j_1+\varepsilon_1}^{-1} x_{i_1+\varepsilon_2}$ cancels out
inside of the product $\delta_{\varepsilon_1}(n)
\delta_{\varepsilon_2}(p)$, and by the inductive assumption we are
done.

If $j_1+\varepsilon_1 < i_1+\varepsilon_2$ then $j_1+\varepsilon_1$
is the smallest index in $\delta_{\varepsilon_1}(n)
\delta_{\varepsilon_2}(p)$ and therefore, using $\CR$, the word
$\delta_{\varepsilon_1}(n) \delta_{\varepsilon_2}(p)$ can be
rewritten the following way:
$$
\begin{array}{c}
\delta_{\varepsilon_1}(n) \delta_{\varepsilon_2}(p) =
x_{j_t+\varepsilon_1}^{-1} \ldots x_{j_2+\varepsilon_1}^{-1}
x_{j_1+\varepsilon_1}^{-1} x_{i_1+\varepsilon_1} \ldots
x_{i_s+\varepsilon_2}
\stackrel{\CR}{\rightarrow} \\
\stackrel{\CR}{\rightarrow} x_{j_t+\varepsilon_1}^{-1} \ldots
x_{j_2+\varepsilon_1}^{-1} x_{i_1+\varepsilon_1+1} \ldots
x_{i_s+\varepsilon_2+1} x_{j_1+\varepsilon_1}^{-1}
\end{array}
$$
Note that since $j_1+\varepsilon_1$ is the smallest index in
$\delta_{\varepsilon_1}(n) \delta_{\varepsilon_2}(p)$, the smallest
index in $w = Merge_{-,+}(x_{j_t}^{-1} \ldots x_{j_2}^{-1},x_{i_2}
\ldots x_{i_s} , \varepsilon_1 , \varepsilon_2)$ is not less than
$j_1+\varepsilon_1$. By the  inductive assumption, $w$ is a
seminormal form for $\delta_{\varepsilon_1}(x_{j_t}^{-1} \ldots
x_{j_2}^{-1}) \delta_{\varepsilon_2}(x_{i_2} \ldots x_{i_s})$.
Therefore, $w x_{j_1+\varepsilon_2+1}^{-1} =_F
\delta_{\varepsilon_1}(n) \delta_{\varepsilon_2}(p)$ and it is a
seminormal form.

The last case where $j_1+\varepsilon_1 > i_1+\varepsilon_2$ is
treated similarly.

\end{proof}

Using ideas from Algorithm \ref{al:semi_norm_form_NP} one can easily
implement an algorithm merging positive words and an algorithm
merging negative words, so that statements similar to Lemma
\ref{le:merging_correct_NP} would hold.
 We will denote these two
algorithms by $Merge_{-,-}(n_1,n_2, \varepsilon_1 , \varepsilon_2 )$
and $Merge_{+,+}(p_1,p_2, \varepsilon_1 , \varepsilon_2 )$,
respectively. Thus, computation of a seminormal form of a product of
two arbitrary seminormal forms has the following form.

\begin{algorithm} \label{al:semi_norm_form}{\em (Seminormal form of a product of seminormal
forms)}\\
{\sc Signature.} $w = Merge(w_1,w_2)$.\\
{\sc Input.} Seminormal forms $w_1$ and $w_2$.\\
{\sc Output.} Seminormal form $w$ such that $w =_F w_1 w_2$.\\
{\sc Computations.}
\begin{enumerate}
 \item[A)] Represent $w_i$ as a product of a positive and negative word ($w_1 = p_1 n_1$ and $w_2 = p_2 n_2$).
 \item[B)] Compute $w' = Merge_{-,+}(n_1,p_2,0,0)$ and represent it as a product of a positive and negative word $w' = p_2' n_1'$.
 \item[C)] Compute $w'' = Merge_{+,+}(p_1,p_2',0,0)$.
 \item[D)] Compute $w''' = Merge_{-,-}(n_1',n_2,0,0)$.
 \item[E)] Output $w'' w'''$.
\end{enumerate}
\end{algorithm}

\begin{lemma} \label{le:merging_correct}
For any pair of seminormal forms $w_1$ and $w_2$ the word $w =
Merge(w_1,w_2)$ is a seminormal form of the  product $w_1 w_2$.
Moreover, the time-complexity of computing $w$ is bounded by
$C(|w_1| + |w_2|)$ for some constant $C$.
\end{lemma}

\begin{proof}
Follows from Lemma \ref{le:merging_correct_NP}.

\end{proof}

\subsection{Seminormal form computation}

\begin{algorithm} \label{al:semi_norm_form2} {\em (Seminormal form)}\\
{\sc Signature.} $u = SemiNormalForm(w)$.\\
{\sc Input.} A word $w$ in generators of $F$.\\
{\sc Output.} A seminormal form $u$ such that $u = w$ in $F$.\\
{\sc Computations.}
\begin{enumerate}
 \item[A)] If $|w| \le 1$ then output $w$.
 \item[B)] Represent $w$ as a product $w_1 w_2$ such that
 $|w_1|-|w_2| \le 1$.
 \item[C)] Recursively compute\\ $u_1 = SemiNormalForm(w_1)$ and\\ $u_2 = SemiNormalForm(w_2)$.
 \item[D)] Let $u = Merge(u_1,u_2)$.
 \item[E)] Output $u$.
\end{enumerate}

\end{algorithm}

\begin{lemma}
Let $w$ be a word in generators of $F$. The output of Algorithm
\ref{al:semi_norm_form2} on $w$ is a  seminormal form for $w$. The
number of operations required for Algorithm \ref{al:semi_norm_form2}
to terminate  is $O(C|w| \log |w|)$, where $C$ is a constant
independent of $w$.
\end{lemma}

\begin{proof}
The first statement can be proved by induction on the length of
$w$. The base of the induction is the  case where $|w| = 1$. In
this case $w$ is already in a seminormal form, and the output is
correct. The induction step was proved in Lemma
\ref{le:merging_correct}.

To prove the second statement we denote by $T(n)$ the number of steps
required for Algorithm \ref{al:semi_norm_form2}  to terminate on an input
of length $n$. Then clearly
$$T(n) = 2T(\frac{n}{2}) + C\cdot n,$$
where the last summand $C\cdot n$ is the complexity of merging two
seminormal forms with the sum of lengths at most $|n|$. It is an
 easy exercise  to show that in this case
$T(n) = O(C \cdot n \log n)$.
\end{proof}

\subsection{Normal form computation}

The next lemma suggests how a pair of generators contradicting the
property (NF2) can be removed and how all  such pairs can be found.

\begin{lemma} \label{le:bad_pair_removal}
Let $w = x_{i_1} \ldots x_{i_s} x_{j_t}^{-1} \ldots x_{j_1}^{-1}$ be
a seminormal form, $(x_{i_a},x_{j_b}^{-1})$ be the pair of
generators in $w$ which contradicts (NF2), where $a$ and $b$ are
maximal with  this property. Let
$$w' = x_{i_1} \ldots x_{i_{a-1}} \delta_{-1}(x_{i_{a+1}} \ldots x_{i_s} x_{j_t}^{-1} \ldots x_{j_{b+1}}^{-1}) x_{j_{b-1}}^{-1} \ldots x_{j_1}^{-1}.$$
Then $w'$ is in a seminormal form and $w =_F w'$. Moreover, if
$(x_{i_c},x_{j_d}^{-1})$ is the pair of generators in $w'$ which
contradicts (NF2) (where $a$ and $b$ are maximal with   this
property), then $c<a$ and $d<b$.
\end{lemma}

\begin{proof}

It follows from the definition of (NF2) and seminormal forms that
all indices in $x_{i_{a+1}} \ldots x_{i_s} x_{j_t}^{-1} \ldots
x_{j_{b+1}}^{-1}$ are greater than $i_a+1$ and, therefore, indices
in $\delta_{-1}(x_{i_{a+1}} \ldots x_{i_s} x_{j_t}^{-1} \ldots
x_{j_{b+1}}^{-1})$ are greater than $i_a$. Now it is clear that $w'$
is a seminormal form. Then doing rewrites opposite to rewrites from
$\CR$ we can get  the word $w'$ from the word $w$. Thus, $w=_F w'$.

There are two possible cases: either $c>a$ and $d>b$ or $c<a$ and
$d<b$. We need to show that the former case is, in fact, impossible.
Assume, by way of   contradiction, that $c>a$ and $d>b$. Now observe
that if $(x_{i_a},x_{j_b}^{-1})$ is a pair of generators in $w$
contradicting (NF2), then
$(x_{i_a+\varepsilon},x_{j_b+\varepsilon}^{-1})$ contradicts (NF2)
in $\delta_{\varepsilon}(w)$. Therefore, inequalities $c>a$ and
$d>b$ contradict  the choice of $a$ and $b$.

\end{proof}

By Lemma \ref{le:bad_pair_removal} we can start looking for bad
pairs in a seminormal form starting from the middle of a word. The
next algorithm implements this idea. The algorithm  is in two
parts. The first part finds all ``bad" pairs starting from the
middle of a  given $w$, and the second part applies
$\delta_\varepsilon$ to segments where it is required. A notable
 feature of Algorithm
\ref{al:semi_norm_form_erase} is that it does not apply the operator
$\delta_{-1}$ immediately (as in $w'$ of Lemma
\ref{le:bad_pair_removal}) when a bad pair is found, but instead, it
keeps the information about how indices must be changed later. This
information is accumulated in two sequences (stacks), one for the
positive subword of $w$, the other one for the negative subword of
$w$. Also, in Algorithm \ref{al:semi_norm_form_erase}, the size of
stack $S_1$ (or $S_2$) equals  the length of an auxiliary word $w_1$
(resp. $w_2$). Therefore, at step B), $x_a$ (resp. $x_b$) is defined
if and only if $\varepsilon_1$ (resp. $\varepsilon_2$) is defined.

\begin{algorithm} \label{al:semi_norm_form_erase}{\em (Erasing bad pairs
from a seminormal form)}\\
{\sc Signature.} $w =EraseBadPairs(u)$.\\
{\sc Input.} A seminormal form $u = x_{i_1} \ldots
x_{i_s} x_{j_t}^{-1} \ldots x_{j_1}^{-1}$.\\
{\sc Output.} A word $w$ which is the  normal form of $u$.\\
{\sc Initialization.} Let $\delta=0$, $\delta_1=0$, $\delta_2=0$,
$w_1 = 1$, and $w_2 = 1$. Let $u_1 = x_{i_1}
\ldots x_{i_s}$ and  $u_2 = x_{j_t}^{-1} \ldots x_{j_1}^{-1}$ be
the positive and negative parts of $u$. Additionally, we set
up two empty stacks $S_1$ and $S_2$.\\
{\sc Computations.}
\begin{enumerate}

\item[A.] Let the  current $u_1 = x_{i_1}
\ldots x_{i_s}$ and $u_2 = x_{j_t}^{-1} \ldots x_{j_1}^{-1}$.

\item[B.] Let $x_a$ be the leftmost letter of $w_1$, $x_b$ the
rightmost letter  of $w_2$, and $\varepsilon_i$ ($i=1,2$) the top
element of $S_i$, i.e., the last element that was put there.
If any of these values does not exist (because, say, $S_i$ is empty),
 then the corresponding variable is not defined.

\begin{enumerate}
 \item[1)] If $s>0$ and ($t=0$ or $i_s>j_t$), then:
\begin{enumerate}
 \item[a)] multiply $w_1$ on the left by $x_{i_s}$ (i.e. $w_1 \leftarrow x_{i_s} w_1$);
 \item[b)] erase $x_{i_s}$ from $u_1$;
 \item[c)] push $0$ into $S_1$;
 \item[d)] goto 5).
\end{enumerate}
 \item[2)] If $t>0$ and ($s=0$ or $j_t>i_s$), then:
\begin{enumerate}
 \item[a)] multiply $w_2$ on the right by $x_{j_t}^{-1}$ (i.e. $w_2 \leftarrow w_2
 x_{j_t}^{-1}$);
 \item[b)] erase $x_{j_t}^{-1}$ from $u_2$;
 \item[c)] push $0$ into $S_2$;
 \item[d)] goto 5).
\end{enumerate}

 \item[3)] If $i_s=j_t$ and (the numbers $a-\varepsilon_1$ and
$b-\varepsilon_2$ (those that are defined) are not equal to $i_s$ or
$i_s+1$), then:
\begin{enumerate}
 \item[a)] erase $x_{i_s}$ from $u_1$;
 \item[b)] erase $x_{j_t}^{-1}$ from $u_2$;
 \item[c)] if $S_1$ is not empty, increase the top element of $S_1$;
 \item[d)] if $S_2$ is not empty, increase the top element of $S_2$;
 \item[e)] goto 5).
\end{enumerate}

 \item[4)] If 1)-3) are not applicable (when $i_s=j_t$ and (one of
the numbers $a-\varepsilon_1$, $b-\varepsilon_2$ is defined and is
equal to either $i_s$ or $i_s+1$)), then:
\begin{enumerate}
 \item[a)] multiply $w_1$ on the left by $x_{i_s}$ (i.e. $w_1 \leftarrow x_{i_s} w_1$);
 \item[b)] multiply $w_2$ on the right by $x_{j_t}^{-1}$ (i.e. $w_2 \leftarrow w_2
 x_{j_t}^{-1}$);
 \item[c)] erase $x_{i_s}$ from $u_1$;
 \item[d)] erase $x_{j_t}^{-1}$ from $u_2$;
 \item[e)] push $0$ into $S_1$;
 \item[f)] push $0$ into $S_2$;
 \item[g)] goto 5).
\end{enumerate}

 \item[5)] If $u_1$ or $u_2$ is not empty then goto 1).
\end{enumerate}

\item[C.] While $w_1$ is not empty:
\begin{enumerate}
 \item[1)] let $x_{i_1}$ be the first letter of $w_1$ (i.e. $w_1 = x_{i_1} \cdot
 w_1'$);
 \item[2)] take (pop) $c$ from the top of $S_1$ and add to  $\delta_1$ (i.e. $\delta_1 \leftarrow \delta_1+c$);
 \item[3)] multiply $u_1$ on the right by $x_{i_1-\delta_1}$  (i.e. $u_1 \leftarrow u_1 x_{i_1-\delta_1}$);
 \item[4)] erase $x_{i_1}$ from $w_1$.
\end{enumerate}

 \item[D.] While $w_2$ is not empty:
\begin{enumerate}
 \item[1)] let $x_{j_1}^{-1}$ be the last letter of $w_2$ (i.e. $w_2 = w_2' \cdot x_{j_1}^{-1}$);
 \item[2)] take (pop) $c$ from the top of $S_2$ and add to $\delta_2$ (i.e. $\delta_2 \leftarrow \delta_2+c$);
 \item[3)] multiply $u_2$ on the left by $x_{j_1-\delta_2}^{-1}$ (i.e. $u_2 \leftarrow x_{j_1-\delta_2}^{-1} u_2$);
 \item[4)] erase $x_{j_1}^{-1}$ from $w_2$.
\end{enumerate}
 \item[E.] Return $u_1 u_2$.
\end{enumerate}

\end{algorithm}

\begin{proposition}
The output of Algorithm \ref{al:semi_norm_form_erase} is the normal
form $w$ of a seminormal form $u$. The number of operations required
for Algorithm \ref{al:semi_norm_form_erase} to terminate is bounded
  by   $D \cdot |u|$, where $D$ is a constant  independent of $u$.
\end{proposition}

\begin{proof}
The first statement follows from Lemma \ref{le:bad_pair_removal}.
The time estimate is obvious from the algorithm since the words
$u_1, u_2$ are processed letter-by-letter, and no letter is processed
more than once.
\end{proof}

 As a corollary, we get the main result of this section:

\begin{theorem}
In Thompson's group $F$, the  normal form  of a given
 word $w$ can be computed in time $O(|w| \log |w|)$.
\end{theorem}


\baselineskip 11 pt


\begin{thebibliography}{ABC}


\bibitem{AAG}
I. Anshel, M. Anshel, D. Goldfeld, {\it  An algebraic method for
public-key cryptography}, Math. Res. Lett. {\bf 6} (1999),
287--291.

\bibitem{CFP}
J. W. Cannon, W. J. Floyd, and W. R. Parry, {\it Introductory notes on Richard
Thompson's groups}, L'Enseignement Mathematique (2) {\bf 42} (1996),
215--256.

\bibitem{CKLHC}
J. C. Cha, K. H. Ko, S. J. Lee, J. W. Han, J. H. Cheon, {\it An
Efficient Implementation of Braid Groups}, ASIACRYPT 2001, Lecture
Notes in Comput. Sci.  {\bf 2248} (2001), 144--156.

\bibitem{HS}
D. Hofheinz and R. Steinwandt, {\it A practical attack on some braid
group based cryptographic primitives}, in Public Key Cryptography,
6th International Workshop on Practice and Theory in Public Key
Cryptography,
 PKC 2003 Proceedings, Y.G. Desmedt, ed., Lecture Notes in Computer Science
{\bf 2567},  pp. 187--198, Springer, 2002.

\bibitem{HT}
J. Hughes and  A.~Tannenbaum, {\it  Length-based attacks for certain
group based encryption rewriting systems}, Workshop SECI02
Securit\`e de la
 Communication sur Intenet, September 2002, Tunis, Tunisia. \\
http://www.storagetek.com/hughes/

\bibitem{KLCHKP}
K. H. Ko, S. J. Lee, J. H. Cheon, J. W. Han, J. Kang, C. Park,
{\it New public-key cryptosystem using braid groups}, Advances in
cryptology---CRYPTO 2000 (Santa Barbara, CA), 166--183, Lecture
Notes in Comput. Sci.  {\bf 1880}, Springer, Berlin, 2000.

\bibitem{Shpilrain}
V. Shpilrain, {\it Assessing security of some
   group based cryptosystems}, Contemp. Math., Amer. Math. Soc.
{\bf 360} (2004), 167--177.

\bibitem{SU}
V. Shpilrain and A. Ushakov, {\it The conjugacy search problem in
public key cryptography: unnecessary and insufficient}, Applicable
Algebra in Engineering, Communication and Computing, to appear.\\
http://eprint.iacr.org/2004/321/

\bibitem{SZ}
V. Shpilrain and G. Zapata, {\it Combinatorial group theory and
public key cryptography}, Applicable Algebra in Engineering,
Communication and Computing, to appear.

\bibitem{Sims}
C. Sims,  {\it  Computation with finitely presented groups},
Encyclopedia of Mathematics and its Applications, {\bf 48}.
Cambridge University Press, Cambridge, 1994.


\end{thebibliography}
\end{document}